# An algorithm for minimization of arbitrary generic functions in one dimension over a finite domain

Glauco Masotti [(*)]

**Abstract**

A new algorithm for one-dimensional minimization is described in detail and the results of some tests on practical cases are reported and illustrated. The method requires only punctual computation of the function, and is suitable to be applied in "difficult" cases, that is when the function is highly irregular and has multiple sub-optimal local minima. The algorithm uses quadratic or cubic interpolation and subdivision of intervals in golden ratio as a last resort. It improves over Brent's method and similar ones in several aspects. It manages multiple local minima, takes into account the complications of having to deal with a finite domain, rather than an unlimited one, and has a slightly faster convergence in most cases.

## 1      One dimensional minimization, improving over classic methods

This procedure for one dimensional minimization is used as a fundamental building block of an optimizer for generic functions in multidimensional space [1]. We developed for this purpose an original algorithm, which improves over Brent and similar methods [2],[3]. In fact multiple local minima are managed, and the complications of having to deal with a finite domain rather than an unlimited one are taken into account, plus a slightly faster convergence has been observed in most cases.
Given *f* the function to be minimized, which should be a real function defined over a domain *[xinf, xsup]*, the algorithm starts taking as input two distinct points *xa* and *xb*, ordered in such a way that going from *xa* to *xb* we go in the downhill direction, i.e. if  *fa=f(xa),* and *fb=f(xb)*, we should have:

$$fb < fa \qquad (1)$$

Two such points can be already available as the output of other processes. If we have no point available we can take the first point *xa* as a random point inside the domain. If we have just one point to start with, the second point *fb* can be computed taking a random step. Obviously if *fb > fa*  it is sufficient to swap the two points to satisfy condition (1). If it happens to be *fa=fb*, we may keep taking other random steps, until we have *fa≠ fb*, or a maximum number of trials is reached, in this case we should assume that the function is a constant.

The algorithm initially activates a procedure which searches in the downhill direction, and in some cases also in the uphill direction, for triplets of points that bracket a local minimum of the function.
Taking as a fundamental reference the popular "Numerical Recipes" book [2], the purpose of this procedure is thus equivalent to the *mnbrak* routine, as described there, however our procedure is capable to deal with a finite domain (which is the cause of some complications), not just an unlimited one, plus we manage the case of functions with multiple local minima, not just one.
Local minima are at first only approximated, then our procedure keeps track of all of them, so that with further processing they are "refined" and pushed closer to the real local minima. At the end the best local minimum will be returned as the global optimum.

## 2      Basic steps of our algorithm

In summary our algorithm can be divided in two phases as follows.

*A) Exploration of domain and construction of a polygonal approximation of the function*

In the 1st phase the domain is explored, taking steps which magnify the initial interval. Assuming that the function can also be erratic (i.e. not smooth), a choice is made in favor of accuracy rather than efficiency, thus steps are normally taken  using a "short" golden ratio magnification, e.g.  if *[a, b]* is the current interval and *[a, c]* is the new one, we take *c=b+0.618(b – a)*, not like in the classic *mnbrak*, which uses a factor of *1.618*.
While decreasing values of the function are found, the exploration continues until the boundary is reached or larger values of the function are found.

If the function is monotonically decreasing the best point (which lays on the domain boundary) is returned.

As the function may have multiple local minima the possibility of having the function decreasing again after rising is contemplated, thus the exploration stops only if a rise is confirmed, i.e. we have two consecutive steps up, otherwise if we have a step down following a step up the exploration continues.

If the initial interval is big with respect to the entire domain, we may have limited space for expanding the initial interval, in this case the initial interval will be a relatively large part of the domain, where the behavior of the function is left unknown. More generally there may be cases where the ratio of the sizes of the biggest interval over the smallest one is large. In these cases an exploration also of the interior of the biggest interval is carried out, using the same criteria of subdivision in golden ratio (note that this would mean moving uphill, given the current knowledge of the function).

All performed calculations of the function are saved for possible future uses in an ordered list, which constitutes a *polygonal* of known values of the function. Unfortunately the behavior of the function in between the calculated points remains unknown, thus if the function is very erratic we can be left unaware of some minima and miss them.

*B) Finding local minima using parabolic and cubic interpolations*

At the end of the exploration process the polygonal is scanned and each triplet which encloses a minimum is processed for "refinement" of this minimum. This is done using parabolic or cubic interpolation in alternative, this latter is attempted only if the former fails (i.e. it does not improve the current minimum). Subdivision of intervals in golden ratios is performed as a last resort (the rationale of using golden ratios can be found in [2]).

The added use of cubic interpolation usually improves performances over the strict use of parabolas, like for Brent's method [2] and similar ones [3].

A comparison of a classic method which uses only parabolas, with our methods which also uses cubic curves, is illustrated in Fig. 1 and Fig. 2. The same portion of a function (represented by the black curve) is considered, this is singled out by the triplet {a, b, c} of the polygonal (shown with a dotted gray line), which form a local minimum in its central point.

The operations of a classic method are illustrated in Fig. 1. The first parabola (shown in magenta) interpolates the three known points {a, b, c} of the function, the abscissa of the minimum is determined and the function is computed at this point, which is indicated as point 1. This is a high value, but it is closer to the current minimum, the second parabola (in orange) interpolates the new triplet {a, b, 1} which identifies this minimum. Point 2 is in correspondence of the minimum of the new parabola and represent a new minimum for the function. The third parabola is illustrated in green and produces point 3. At this point we should interpolate the points {a, 3, 2}, but it happens that this parabola has a minimum practically coincident with point 3, so that it is not possible to compute a new parabola. Point 4 is thus produced by subdividing in golden ratio the interval [3, 2] of the two points lower in value. By interpolation of points {3, 4, 2} we obtain point 5, then interpolating {3, 5, 4} we have no new point again, thus point 6 is computed by subdividing the interval [3, 5]. Points 6 and 5 are sufficiently close in abscissa and ordinate to call for convergence achieved and determine the end of the algorithm. Therefore with this method, in this example, we had to compute the function six times in order to find the minimum.

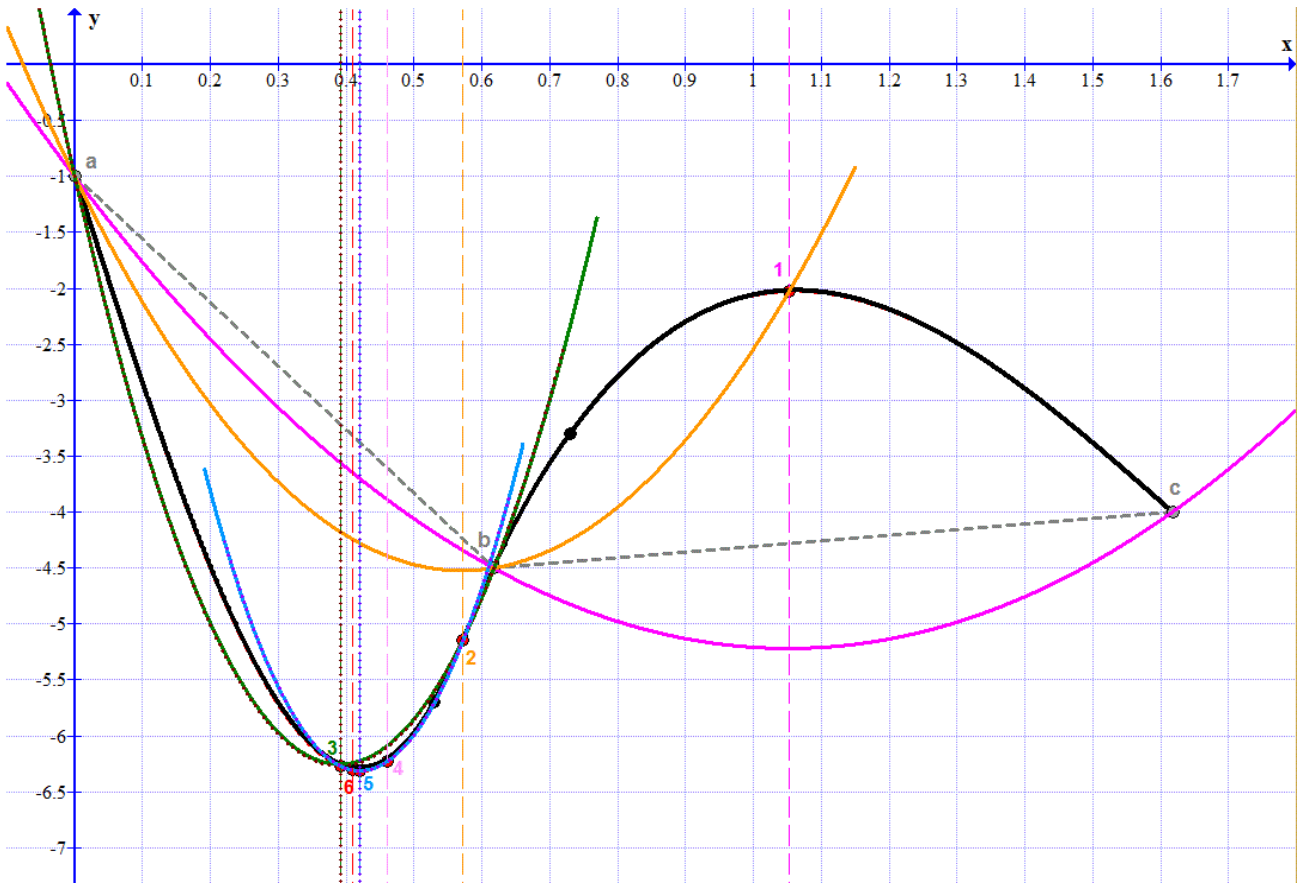

**Figure 1.** Finding a minimum with successive parabolic interpolations.

The method which was developed and implemented in our algorithm, applied to this example case, is illustrated in Fig. 2 and proceeds as follows.
The first step is identical to the classic method, thus the parabola (shown in magenta) which interpolates the three known points {a, b, c} and the abscissa of its minimum are determined and the function is computed at this point, which is shown as point 1.
This point does not improve over the current best point (the minimum at b), therefore the second step uses a cubic, interpolating the points {a, b, 1, c}. The parabolic interpolation would have used the three points closer to the current minimum, while the cubic interpolation requires four points, thus also point c is used. The cubic conveys more information on the behavior of the function, in fact it approaches the function more closely throughout the interval [a, c], although we are interested only in a better approximation in the subinterval [a, 1] (where the minimum is located), but there are chances to get also this as a byproduct. The cubic interpolation leads to point 2, which improves over b.
The third step uses a parabola (which is always attempted as a first trial with new points) interpolating points {1, 2, b}, in correspondence of the abscissa of its minimum point 3 is found, this does not improve over point 2, thus a cubic through {1, 3, 2, b} is attempted, finding point 4, which improves over the minimum of point 2. Points 4 and 2 are close enough to each other in abscissa and ordinate values, so that convergence is achieved and the algorithm is terminated.
In this case our algorithm requires only four iterations, with respect to the six iterations required by the classic method, certainly this is not enough to assess the superiority of our mixed parabolic/cubic approach, however the results of our tests seem to show that our approach is statistically convenient.

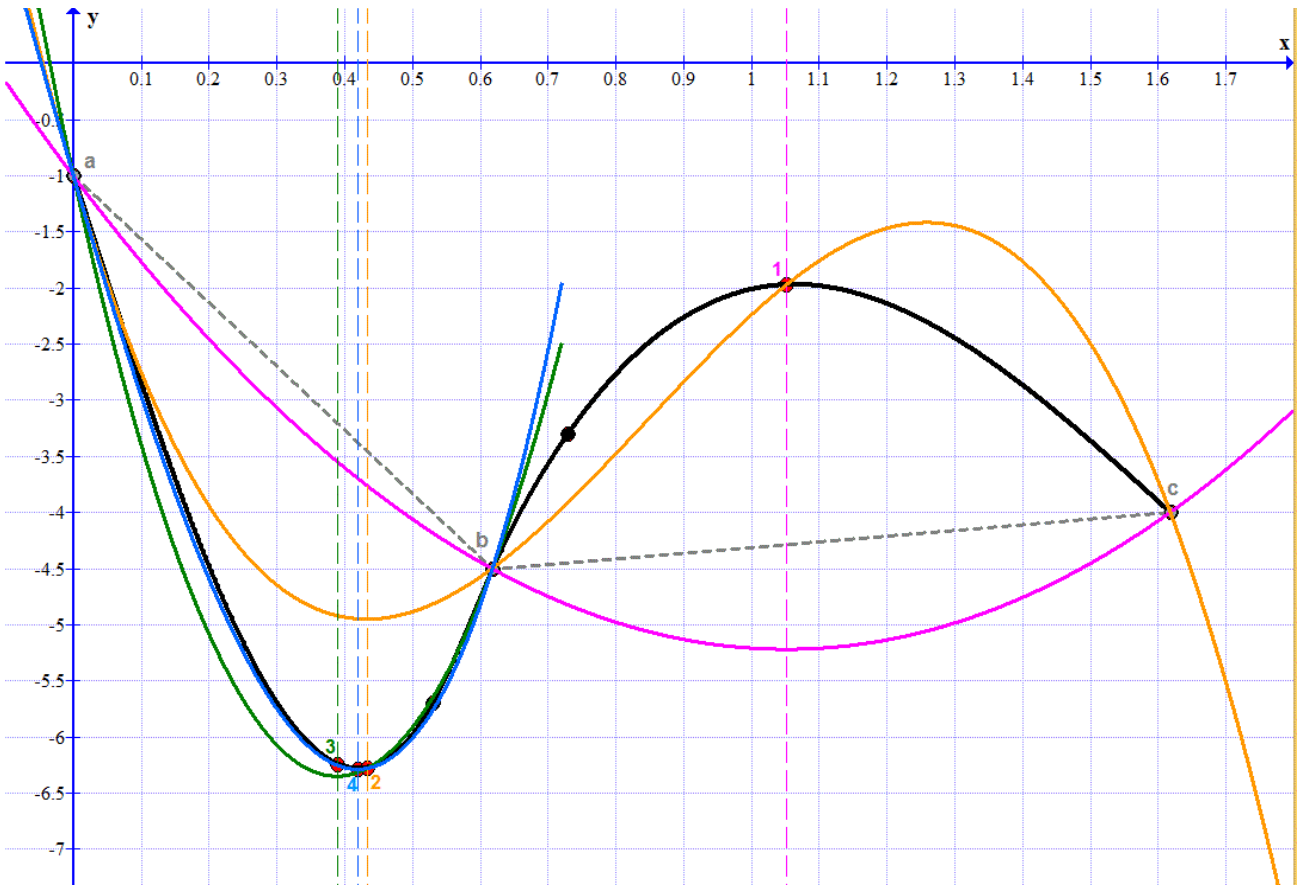

**Figure 2.** Finding a minimum using parabolic and cubic interpolations.

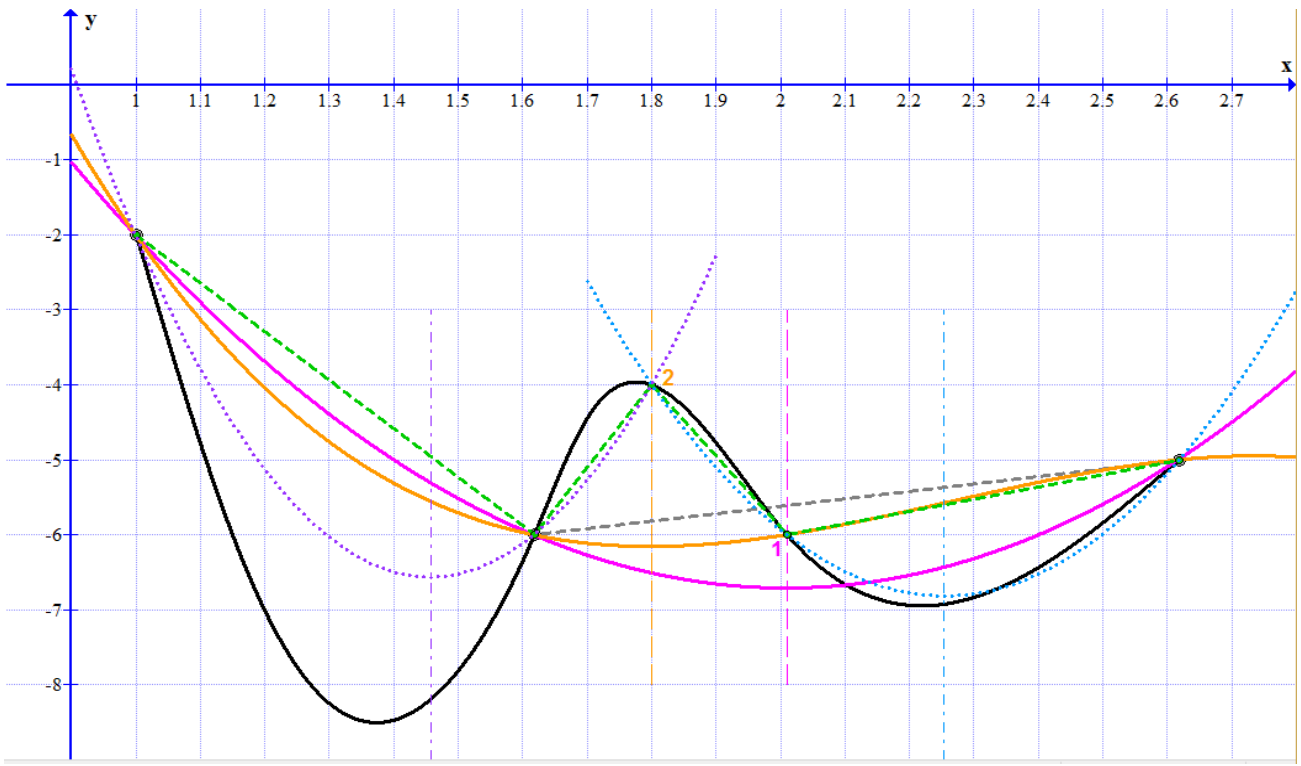

**Figure 3.** The refinement of the first local minimum reveals the existence of another local minimum.

In our algorithm nothing is wasted, the process of minimum refinement adds points to the polygonal, which may then indicate new local minima of the function, which need further refinements. Figure 3 shows a case of this kind. The initial polyline is shown in gray, after finding point 1 and 2 by means of parabolic and cubic interpolation, the updated polygonal, shown in green, has two local minima, which need separate refinement. The two initial parabolic interpolations for each minimum are shown with dotted lines.

The algorithm, in our last implementation, uses only the two phases of above, but we also tested a variant. The rationale of this variant is as follows.
Is there anything better that we can do, other than expanding the initial interval, adding intervals in golden ratio, like in phase one? The *mnbrak* routine, as described in [2], also uses parabolic extrapolation for predicting the behavior of the function and thus evaluates the function also at the minima of the parabolic extrapolations. Unfortunately parabolic extrapolation should seldom represent well the function, like in most of the cases that we have treated, thus most of the times it is useless and we preferred not to attempt it.
However, the polygonal which comes out from phase one can be too coarse, especially at the beginning, thus a refinement of it in smaller intervals may be advisable, in fact a coarse polygonal may not reveal the existence of potential minima. An example of this is given in Fig. 4.

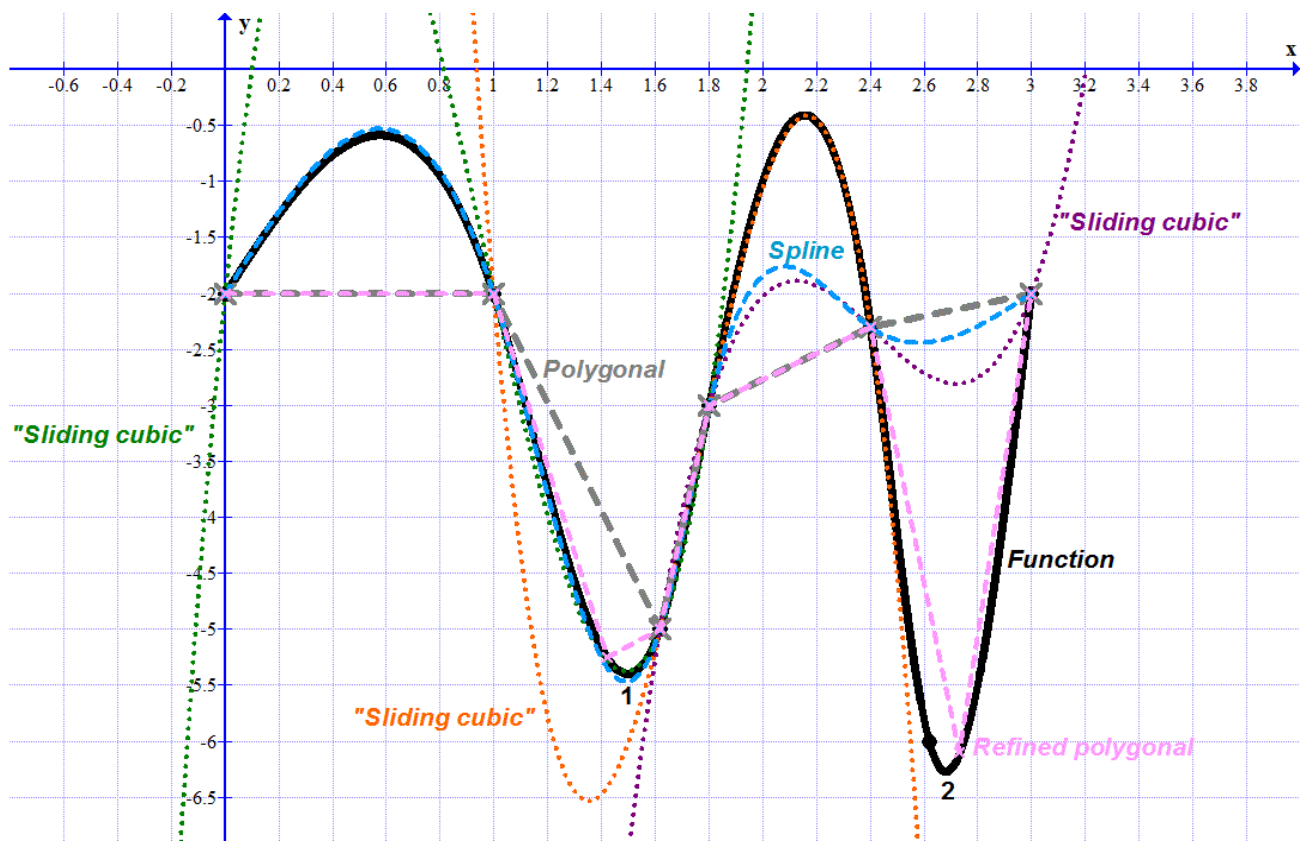

**Figure 4.** Polygonal refinement as an optional stage.

A function is shown, represented by a black solid line, and an approximating polygonal with five sides is shown with a dashed gray line. We see that this polygonal identifies a subinterval which contains the local minimum 1, but does not reveal the existence of the global minimum in 2.
To detect this situations a spline interpolating all points could be used (this is shown with a dashed light blue line). In fact, although the function has much wider oscillations, also the spline has a local minimum in the same subinterval of 2. Rather than using a spline we implemented instead a simpler technique, which most of the times gives qualitatively equivalent results. We used what we called a "sliding" cubic interpolation, which is a local interpolation with a cubic of each subinterval of four points. These cubics are shown in dotted lines of different colors. The suspect local minimum locations are calculated and the function is checked at these points, which are added to the polygonal. The refined polygonal is shown with a dashed pink line. In this case the refined polygonal is able to identify both intervals of potential minimum and the effort was worth it, but we have to point out that, if the function where much smoother, there would not be any minimum other than 1, thus the predictions would be false and this extra work useless. But also if the function were much more erratic the predictions using the cubic interpolations could not be verified.

In practice this procedure was not much convenient in our test cases, thus it was opted out in the final version of our algorithm.

## 3   Examination of some real cases

For debugging purposes we provided the implementation of our algorithm with extensive graphic feedback, thus we are able to show some real cases which we encountered in using the algorithm for line minimization in the context of optimization of complex functions of several variables.

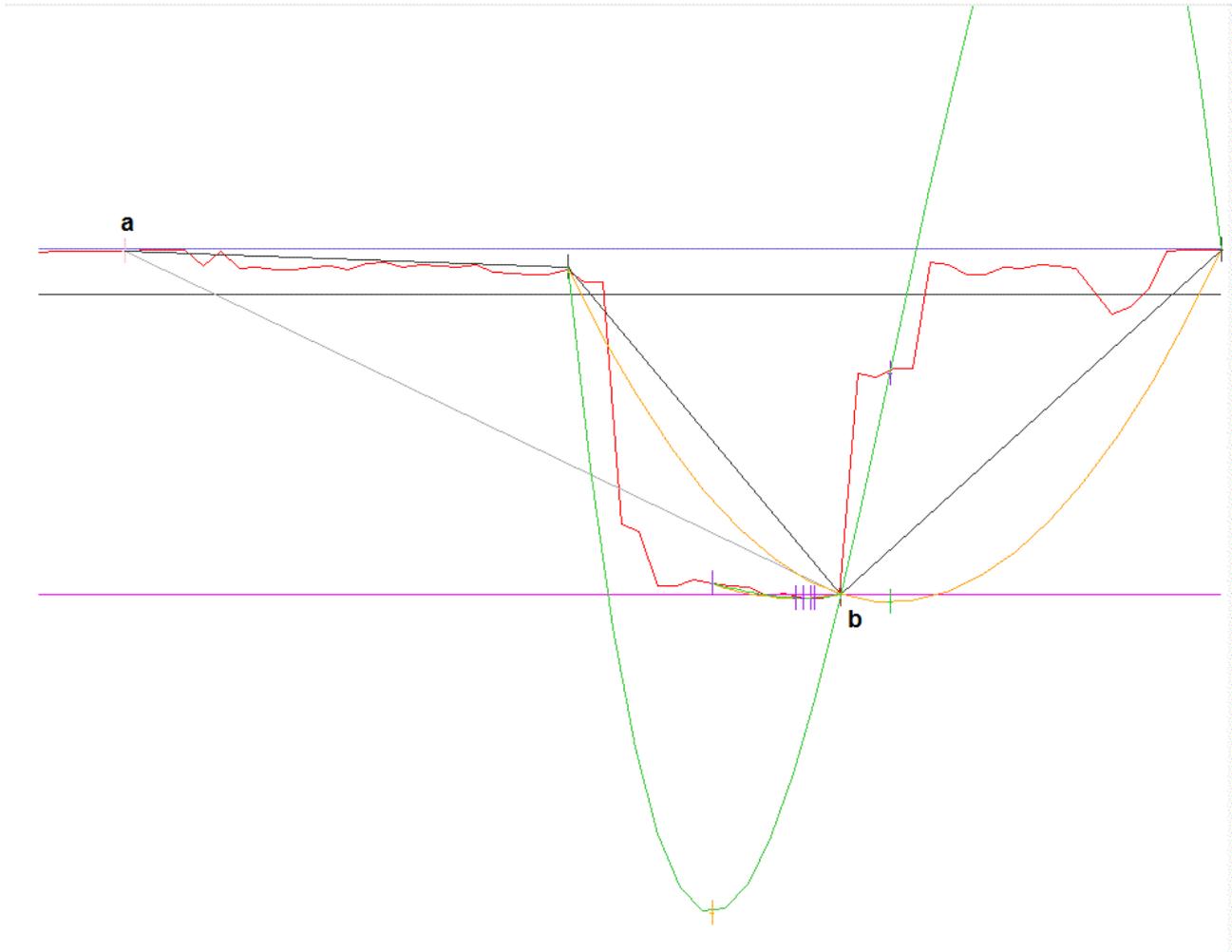

**Figure 5.** A simple real case.

In Figure 5, a simple real case is shown. The function has been sampled in several tens of points, the corresponding polygonal approximation is shown in red. The initial segment, joining the two initial points *{a, b}*, is shown in light gray. The coarse polyline, made up of only three segments, resulting after the first phase of the algorithm is shown in dark gray. The computed interpolating parabolas are shown in orange, while their minimum points are marked in green. These color are swapped for the cubic interpolants. The points of the function which are computed by the algorithm are shown in violet. We see that the global minimum is precisely found although, unfortunately in this case, it is not much better than the initial point *b*.

A somewhat more complex case is shown in Figure 6. The same colors for the various elements are used.

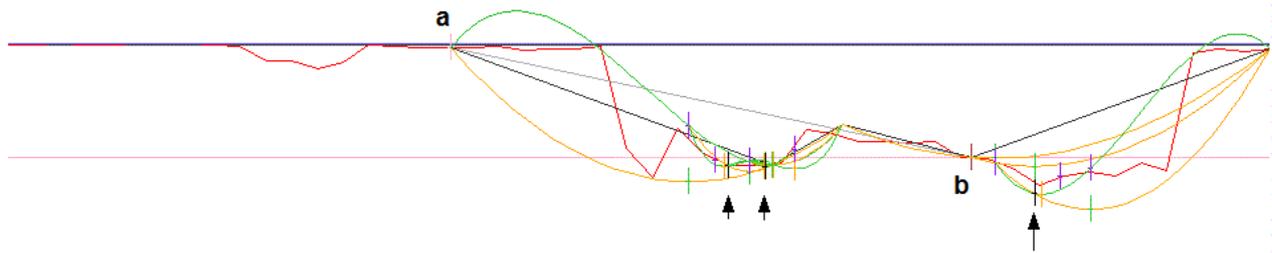

**Figure 6.** All steps of the algorithm for this case are resumed. Three local minima are found.

Here the initial polyline has two local minima. Refining them a third local minimum is discovered. These local minima are marked in black and pointed by the black arrows. Obviously the lowest one is returned.

Figure 7 shows a case where the function in nearly constant for the lower part of the domain, then it starts oscillating wildly approaching the upper end. The polyline connecting the sampled points of the function is shown in red, the initial segment *{a, b}* is shown in light gray, the final polyline connecting all evaluated points of the function (marked with a black cross) is shown in black. The three black arrows point to the three local minima found, two of which improve considerably the initial minimum at b.

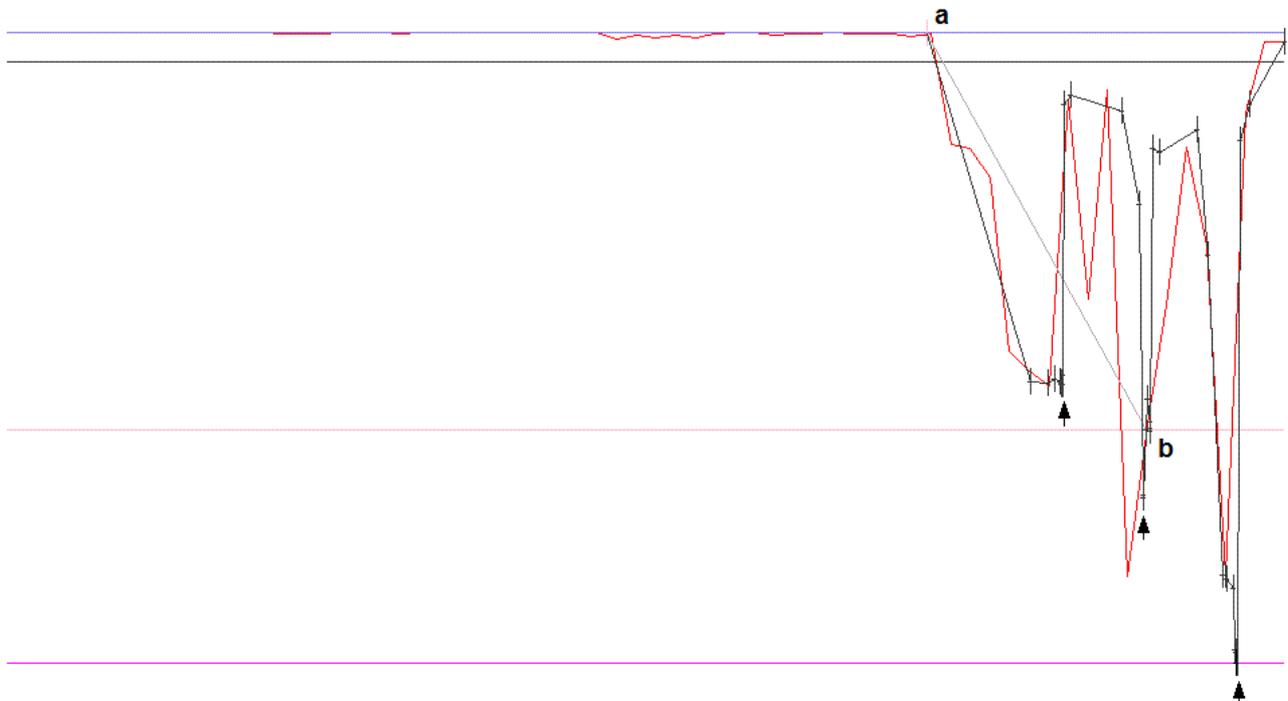

**Figure 7.** Here the function oscillates wildly, but the global minimum is found.

## 4  Using bounds

It may be convenient to use bounds to avoid the refinement of minima which are not promising enough. We should be cautious in doing this, because on one side this may save many unuseful computations of the function most of the times, but, on the other side, it can cause sometimes missing some good minima. An example were the use of bounds is convenient is shown in Figure 8.
A simple bound that can be used is a threshold on the current value of the local minimum of a polygonal line: if the value is below the threshold it is refined (i.e. we search for a better approximation of the real minimum) otherwise we do nothing. It is reasonable to make the bound dynamic, that is while the algorithm proceeds, and the current global minimum is lowered, also the bound is lowered in some way.
In the case of Figure 8 we have initially a polygonal with four local minima below the current bound (the threshold) which is represented by the horizontal dark gray line, thus they are all good candidates for the global minimum. But the refinement of minimum 1 leads to point 1', which is much lower. The updated bound is shown in green, it happens that all other local minima 2, 3, and 4 are above the updated bound so that they are not refined and this saves a lot of unuseful work, because refinement of these points would not improve the current optimum.

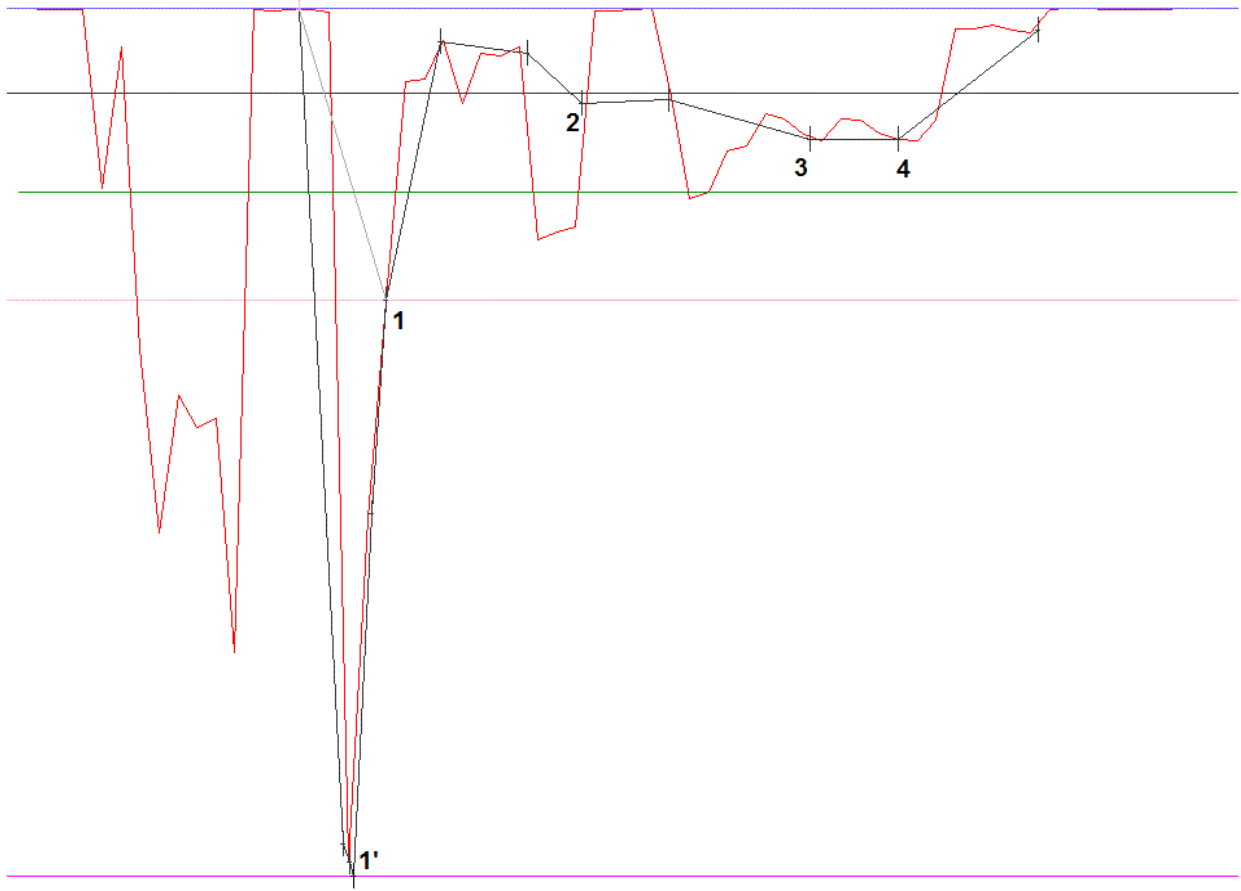

**Figure 8.** In this case it is convenient to use bounds to save unuseful computations.

We found convenient to use other bounds, one is on the absolute variation expressed by the triple $\{p_{k-1}, p_k, p_{k+1}\}$ of a polygonal that identify a candidate point $p_k$ for refinement and another on the slope of the two segments. This means that the finite difference and the difference quotient of the triple should be enough "promising" in order to proceed with refinement, otherwise the chances are high that the function is quasi-constant, so that it would be a waste of resources to compute the refinement of $p_k$ until convergence is achieved. Another bound that we have enforced is on the maximum number of function evaluations done without improving the minimum, exceeding the prefixed limit no other attempts are made and the best point found so far is returned.

## 5    The algorithm

Rather than presenting a synthetic description of the algorithm in some sort of pseudocode, we decided to publish a great part of the original code of the algorithm, which has been implemented in plain C. For sake of clarity, not to disperse attention in too many details, and also not to use too much space, not every detail of the code is reported, but we think that more than enough details are given to fully understand the algorithm, so that, if required, a working implementation should be arranged without too much effort. In particular the code for the management of a list of structures is not reported in detail. The procedures required can be implemented in different, but functionally equivalent, ways. For efficiency reasons, in terms of time and space, we preferred an implementation in the form of a dynamic array of structures, but also an implementation as a linked list of structures may work.
It may be worth noting that we used the Boehm's garbage collector [4] for memory management, which, like in other occasions, simplified somewhat our work, however the algorithm is not too much demanding in terms of memory management, so that also a classic scheme, based on malloc and free, is suitable.

Also note that assertions are commented out, but they are valid, so that, if the program works correctly, they must be verified.

The algorithm uses some numeric tolerances in various parts, that should be tailored for the problem at hand. A value of $10^{-6}$ for the relative tolerance should work in most cases. In general a relative tolerance should be used for the variable and for the function, so that, regardless the range of values assumed by these entities, the results can be calculated with a certain number of meaningful significant digits, but in case the domain of the parameter is normalized an absolute tolerance may be used for it as well.

An approach of this kind can be satisfying for most applications. Relative tolerances should account for the intrinsic "granularity" of the problem at hand (i.e. it makes no sense calculating results to a precision much greater than what can be appreciated in practice), but also for the accumulation of round-off errors. In well-conditioned problems the round-off errors are orders of magnitude smaller than the granularity of the problem, so that a simple heuristic technique of this kind is suitable, otherwise a more sophisticated approach is necessary [5].

```c
// Constant and macro definition:

#define   GOLD 0.61803398875
#define  CGOLD 0.38196601125
#define   EPS2 1e-12
#define   HUGE 2e20
#define TO_BE_EVALUATED (HUGE/2.)
      /* A conventional and unlikely value, which signals that the function
         has not yet been evaluated for the corresponding value of the parameter */

#define SQR(X) ((X)*(X))
#define SQR_(T, X) (T = (X), T*T)

#define   STRUCT void
#define  PSTRUCT void *

// Definition of special structures:

typedef struct s_fPoint {
// An element of the array which defines the polygonal of known points of the function
   double x;  double y; // abscissa and ordinate of this point
   short int marked; /*  1 : the point has already been "refined"
                         0 : otherwise */
} FPOINT, *PFPOINT;

typedef struct s_structList {
// Header of dynamic array of structures,
// used to represent the polygonal approximation of the function
   int step;    // step of allocation, in number of elements
   int size;    // allocated area, in number of elements
   int len;     // number of used (active) elements
   int cursor;  // index of the current element
   void *root;  // pointer to the array of structures
} STRUCTLIST, *PSTRUCTLIST;

// Global variables:

static int nff; // total number of functions evaluations so far
static double xinf, xsup; // limits of the domain of the function

int dontWasteTimeIn_minSearch = 1; // if set enables the use of bounds
double deltaBound; /* minimum variation required on the guess points
                      in order to proceed with the algorithm (set dynamically) */
double slopeBound; /* minimum slope required on the guess points
                      in order to proceed with the algorithm (set dynamically) */
int nMaxFailed = 4; // maximum number of function evaluations without an improvement
double k_ysup = 0.5; /* used to define the upper bound
                        above which local minima are not refined */

// Functions and procedures used by the algorithm:

static double(*ff)(double); // The function!
static double f1(double x); // calls ff and updates the number of evaluations

int findMin1d // subroutine (see below)
```

```c
                 (double xa, double xb, double fa, double fb, double xtol, double ftol,
                  double *xmin, double *ymin);
int quadCubeMinSearch // subroutine (see below)
                    (FPOINT *pA, FPOINT *pB, FPOINT *pC,
                     STRUCTLIST *poly, STRUCTLIST *poly1, double xtol, double ftol,
                  double *xlocmin, double *ylocmin, double *xmin, double *ymin);
int polyAddPoint // subroutine (see below)
               (STRUCTLIST *poly, double x, double *y, int m, double xtol);
double adjustToBoundary // function (see below)
                       (double xa, double xb, double xc, double minRatio, int *clip);
int reqr(double a, double b, double eps); // function (see below)

void structListCreate(PSTRUCTLIST header, unsigned int szel);
    /* Creates a dynamic array of structures of "szel" elements,
       the structure "header" is filled and returned. */
void structListKill(PSTRUCTLIST header);
    /* Deletes the list, no more used, pointed by "header" */
void structListInsertElement(PSTRUCTLIST header, PSTRUCT elem, unsigned int szel);
    /* Insert the element "elem", in the list pointed by "header",
       prior to the element pointed by the cursor */
int  structListRead(PSTRUCTLIST header, PSTRUCT elem, unsigned int szel, int incr);
    /* Places in "elem" the element of the list pointed by "header"
       and addressed by the index cursor+incr */
void structListWrite(PSTRUCTLIST header, PSTRUCT elem, unsigned int szel);
    /* Overwrite the element pointed by the cursor with the contents of "elem" */
int  structListReadNext(PSTRUCTLIST header, PSTRUCT elem, unsigned int szel);
    /* Increments the cursor and returns in "elem" the element pointed by this. */
void structListCopy(PSTRUCTLIST header, PSTRUCTLIST hcopy, unsigned int szel);
    /* Copies in "hcopy" the list pointed by "header" */

//-------------------------------------------------------------------------
// Main procedure of the algorithm:

int
minSearch1d(
    double xa, double xb, double fa, double fb, // initial points
    double(*f)(double), // function to be optimized
    double xtol, /* minimum significant relative variation of the variable
                    (i.e. approximation required on the parameter) */
    double ftol, /* minimum significant relative variation of the function
                    (i.e. approximation required on the minimum value of the function) */
    double x_inf, double x_sup, // domain of the variable (aka parameter) of the function
    double *xmin, // initial guess for the minimum and final value found
    double *ymin  // corresponding value of the function
)
/* Searches for a minimum of f, defined over the domain [x_inf, x_sup],
in a neighborhood of the two initial data points {xa, xb}, which correspond to the values of the function fa
and fb. It is assumed that fb <= fa.
If "bounds" are enforced, it is also required at least a minimal significant variation
between fa and fb (otherwise the procedure returns immediately).
Returns the number of function evaluations */
{
    int nf = 0;

    //ASSERT(xsup > xinf);
    //ASSERT(contained(xa, xinf, xsup) && contained(xb, xinf, xsup));
    //ASSERT(fa < TO_BE_EVALUATED);
    //ASSERT(fb < TO_BE_EVALUATED);
    //ASSERT(fb <= fa);
    //ASSERT(*ymin < TO_BE_EVALUATED);
    //ASSERT(*ymin == fb && *xmin == xb);
    //ASSERT(deltaBound < 0.);
    //ASSERT(slopeBound < 0.);

    xinf = x_inf, xsup = x_sup;
    ff = f;
    if (dontWasteTimeIn_minSearch // bounds are enforced
        && fb - fa > deltaBound /* fb ~ fa */
        || (fb - fa) / (fabs(xb - xa) + EPS2) > slopeBound)
        // The function does not show enough variation, probably better to give up!
```

```
        goto end;
    nf = findMin1d(xa, xb, fa, fb, xtol, ftol, xmin, ymin);

end:
    return nf;
}

//----------------------------------------------------------------------
// Most of the work is done here:

int
findMin1d(double xa, double xb, // extremes of the initial interval
    double fa, double fb, // corresponding values of the funtion, with fb < fa
    double xtol, // tolerance (relative precision) to be used for x,
                 // i.e. minimum variation in x for two points to be considered as distinct
    double ftol, // tolerance (relative precision) to be used for y
    double *xmin, double *ymin // current minimum and corresponding value of the function
)
/*
Given a function f, defined over a domain [xinf, xsup] (static globals),
and given distinct initial points (xa, fa) and (xb, fb), with fb < fa,
this routine searches, in the downhill direction, and in case also in the uphill direction,
for triplets of points that bracket local minima of the function.
The best minimum found is placed in (xmin, ymin)
The number of function evaluations is returned
*/
{
    // the initial interval is saved for possible later use
    double xaa = xa, faa = fa;
    double xbb = xb, fbb = fb;
    double xc, fc, xd, fd;
    int clip, changes;
    int nvertex = 0;
    int expansionCount = 0; // counts the number of expansions of the initial interval
    double const minRatio = SQR(CGOLD); // = 0.145898 used in adjustToBoundary
    double K = GOLD; /* (1+K) is the amplification factor of the initial interval,
                       used to determine the III point */
    FPOINT p1, p2, p3; // used to store a triplet of points
    STRUCTLIST fPoly; // the polygonal of known points that approximates the function

    nff = 0; // reset the counter of function evaluations

    //ASSERT(fb < fa); ASSERT(*xmin == xb && *ymin == fb);

    // initialize the polygonal with the two initial points
    structListCreate(&fPoly, sizeof(FPOINT));
    nvertex += polyAddPoint(&fPoly, xa, &fa, 0, xtol);
    nvertex += polyAddPoint(&fPoly, xb, &fb, 0, xtol);
    //ASSERT(nvertex == 2); // the two initial points must be distinct!

    do { /* extend the initial interval until we have a confirmed rise of the function
            or the boundary is reached */
        xc = xb + K*(xb - xa); /* (1+K) is the amplification factor of [xa, xb]
                                  used to determine xc */
        xc = adjustToBoundary(xa, xb, xc, minRatio, &clip); /* if xc is close or passed the
boundary it is adjusted or clipped to the boundary limit */
        if (!clip) expansionCount++; // only full expansions are counted
        fc = TO_BE_EVALUATED;
        if (polyAddPoint(&fPoly, xc, &fc, 0, xtol)) {
            nvertex++;
            if (fc < *ymin)
                *xmin = xc, *ymin = fc;
        }
        else { /* the last point has not been added to the polygonal
                  => we reached the domain boundary
                  N.B. Also at the first step we may end up here:
                  if xb is very close to the boundary of the domain or |xb-xa| ~ xtol */
            //ASSERT(reqr(xc, xa, xtol) || reqr(xc, xb, xtol) || xc == xinf || xc == xsup);
            break;
        }

        if (fc >= fb) {
            //ASSERT(fb < fa && fb <= fc);
            // {xa, xb, xc} is a triplet which encloses a local minimum
```

```
            xd = xc + K*(xc - xa); // let's try to expand further
            xd = adjustToBoundary(xa, xc, xd, minRatio, &clip);
            if (!clip) expansionCount++;
            fd = TO_BE_EVALUATED;
            if (polyAddPoint(&fPoly, xd, &fd, 0, xtol)) {
                nvertex++;
                if (fd < *ymin)
                    *xmin = xd, *ymin = fd;
            }
            else {
                //ASSERT(xd == xinf || xd == xsup); // once reached the domain boundary we stop
                break;
            }

            if (fd >= fc) /* the function rose again,
                            no other local minima in sight, better not to insist! */
                break;
            else { // fd < fc, the function is again in descent, better to look further
                xa = xc, fa = fc; // a <-- c
                xb = xd, fb = fd; // b <-- d, restart with a new cycle of expansion
            }
        }
        else { // fc < fb, continue to search downhill
            xb = xc, fb = fc; // b <-- c, restart with a new cycle of expansion
        }
    } while (!clip); // == while (inside(xb, xinf, xsup));

    if (expansionCount < 2) {
        /* The initial segment is big with respect to the domain => it was not possible to
           expand it enough, so that it remains the largest part of the polygonal
           => better to subdivide it and look inside */
        xc = xaa + K*(xbb - xaa); // intermediate point between a and b (closer to b)
        fc = TO_BE_EVALUATED;
        if (polyAddPoint(&fPoly, xc, &fc, 0, xtol)) {
            nvertex++;
            if (fc < *ymin)
                *xmin = xc, *ymin = fc;
        }
        //else noop();
    }

    // Let's take the local minima of the polygonal and "refine" them
    double deltap, deltaq;
    double ytol;
    STRUCTLIST fPoly1; // new polygonal, with the newest points added as we find them
    double xlocmin, ylocmin; // current local minimum of the polygonal
    double fPmin = HUGE, fPmax = -HUGE; // current min e max values of the polygonal

    structListCreate(&fPoly1, sizeof(FPOINT));
    structListCopy(&fPoly, &fPoly1, sizeof(FPOINT)); // fPoly --> fPoly1

    do {
        if (dontWasteTimeIn_minSearch) { // let's find out min and max values first
            fPoly.cursor = -1;
            while (structListReadNext(&fPoly, &p1, sizeof(FPOINT))) {
                if (p1.y < fPmin) fPmin = p1.y;
                if (p1.y > fPmax) fPmax = p1.y;
            }
            //ASSERT(fPmin == *ymin);
        }
        changes = 0;
        fPoly.cursor = -1;
        structListReadNext(&fPoly, &p1, sizeof(FPOINT));
        structListReadNext(&fPoly, &p2, sizeof(FPOINT));
        if (p1.y < p2.y)  xlocmin = p1.x, ylocmin = p1.y;
        else              xlocmin = p2.x, ylocmin = p2.y;
        deltap = p2.y - p1.y;
        while (structListReadNext(&fPoly, &p3, sizeof(FPOINT))) {
            if (p3.y < ylocmin) {
                xlocmin = p3.x, ylocmin = p3.y; // monotone decreasing sequence
            }
            deltaq = p3.y - p2.y; ytol = ftol*(1. + fabs(p2.y));
            if (!p2.marked && deltap*deltaq <= 0.
                && ((!dontWasteTimeIn_minSearch && (deltap < -ytol || deltaq > ytol))
```

```c
                    || (dontWasteTimeIn_minSearch
                        && (deltap < deltaBound || deltaq > -deltaBound)
                        && (deltap / (p2.x - p1.x) < slopeBound
                            || deltaq / (p3.x - p2.x) < slopeBound)
                        && p2.y < fPmax - (fPmax - fPmin)*k_ysup
                        // p2 is in the "lower band" of the values found so far
                        )
                    )
                ) {
                /* p2 is a local minimum of the polygonal which has to be refined
                   => we look for a better approximation of the true minimum in its vicinity   */
                changes = quadCubeMinSearch(&p1, &p2, &p3, &fPoly, &fPoly1, xtol, ftol,
                                            &xlocmin, &ylocmin, xmin, ymin);
                            /* searches for the true minimum in the interval [p1,p3]
                            N.B. The new points found are added to fPoly1 */
            }
            memcpy(&p1, &p2, sizeof(FPOINT));
            memcpy(&p2, &p3, sizeof(FPOINT));
            deltap = deltaq;
        }
        structListKill(&fPoly);
        structListCopy(&fPoly1, &fPoly, sizeof(FPOINT)); // fPoly1 --> fPoly
     // the next cycle is executed over the new polygonal (with the added new points)
    } while (changes);
    // when the polygonal doesn't change anymore we are done!

    return nff; // the total number of function evaluations is returned
}
//-------------------------------------------------------------------------

double quadricMin(FPOINT *p1, FPOINT *p2, FPOINT *p3, double *P2ymin);
double cubicMin(FPOINT *p1, FPOINT *p2, FPOINT *p3, FPOINT *p4, double *P3ymin);

//-------------------------------------------------------------------------
// Refinement of a local minimum:

int
quadCubeMinSearch(FPOINT *pA, FPOINT *pB, FPOINT *pC, // initial triple {pA, pB, pC}
                  STRUCTLIST *poly, STRUCTLIST *poly1,
                  double xtol, double ftol,
                  double *xlocmin, double *ylocmin, // local minimum found
                  double *xmin, double *ymin // current global minimum
)
{   /* searches for the true minimum in the interval from p1 to p3
       (with p2 being the current best known point in the interval),
       returns the number of points added to the polygonal poly1,
       marks the points that converge on the minimum found */
    double xtry, ycast, ytry;
    double ytol;
    int newp, changes = 0;
    int nFailed = 0;
    FPOINT p0, p1, p2, p3, p4;

    memcpy(&p1, pA, sizeof(FPOINT));
    memcpy(&p2, pB, sizeof(FPOINT));
    memcpy(&p3, pC, sizeof(FPOINT));

    *xlocmin = p2.x, *ylocmin = p2.y; // the input point p2 is the initial local minimum

#define READ_TRIPLET(PSL, A, B, C, SIZE) \
    structListRead(PSL, &A, SIZE, -1); \
    structListRead(PSL, &B, SIZE,  0); \
    structListRead(PSL, &C, SIZE,  1);
    /* this reads A, B, C from the list of structures (of dimension SIZE) pointed by PSL,
       with B being the element pointed by the cursor,
       A the preceding element and C the subsequent, the cursor is left unchanged! */

    do {
        ytol = ftol*(1. + fabs(p2.y));
        if (p1.y - p2.y < ytol && p3.y - p2.y < ytol)
            break; /* convergence in y is verified, we are done,
                      if the variation of the function over the three points
                      is below the tolerance level it is useless to insist further */
```

```c
        //ASSERT(p2.y < p1.y && p2.y <= p3.y || p2.y < p3.y && p2.y <= p1.y);
      // the minimum point must always be p2

      xtry = quadricMin(&p1, &p2, &p3, &ycast);
      // computes the minimum of the parabola passing through p1, p2, p3
      //ASSERT(inside(xtry, p1.x, p3.x));
      // the parabola must have a minimum inside the interval p1--p3
      ytry = TO_BE_EVALUATED;
      newp = polyAddPoint(poly1, xtry, &ytry, 0, xtol);
      if (newp) { // the computed point is new => it has been added to the polygonal
         changes++;
         if (ytry >= *ylocmin) {
            /* quadratic interpolation did not produce a better point => we try with
               cubic interpolation, taking into account the last point found */
            nFailed++;
            structListRead(poly1, &p4, sizeof(FPOINT), 0);
            // N.B. the cursor of poly1 still points to the new point just added
            //ASSERT(*xlocmin == p2.x); // p2 is still the minimum
            xtry = p4.x < p2.x ? cubicMin(&p1, &p4, &p2, &p3, &ycast)
                               : cubicMin(&p1, &p2, &p4, &p3, &ycast);
            if (xtry < HUGE) {
               //ASSERT(inside(xtry, p1.x, p3.x));
               ytry = TO_BE_EVALUATED;
               newp = polyAddPoint(poly1, xtry, &ytry, 0, xtol);
               if (newp) {
                  changes++;
                  if (ytry < *ylocmin) {
                     // OK, cubic interpolation has found a new local minimum!
                     *xlocmin = xtry, *ylocmin = ytry;
                     if (ytry < *ymin) { // the local minimum is also the new global minimum!
                        *xmin = xtry, *ymin = ytry;
                        nFailed = 0; // reset the counter of failed attempts
                     }
                     else
                        nFailed++; // local minimum has not been improved

                     // cursor still points to the new point just added
                     READ_TRIPLET(poly1, p1, p2, p3, sizeof(FPOINT)); // updates p1, p2, p3
                     p2.marked = 1; structListWrite(poly1, &p2, sizeof(FPOINT));
                     // current local minimum is marked
                     continue; // proceed with a new cycle
                  }
               }
            }
            else /* cubic has not a minimum in the useful domain
                    => generate the new point with subdivision in golden ratio */
               goto gold1;
         }
         else { // OK, quadratic interpolation has found a new local minimum!
            *xlocmin = xtry, *ylocmin = ytry;
            if (ytry < *ymin) { // the local minimum is also the new global minimum!
               *xmin = xtry, *ymin = ytry;
               nFailed = 0; // reset the counter of failed attempts
            }
            else
               nFailed++; // local minimum has not been improved

            READ_TRIPLET(poly1, p1, p2, p3, sizeof(FPOINT)); // updates p1, p2, p3
            p2.marked = 1; structListWrite(poly1, &p2, sizeof(FPOINT));
            // current local minimum is marked
            continue; // proceed with a new cycle
         }
gold1:
         /* if we are here quadratic and cubic interpolation have failed
            => we try the subdivision of the interval in golden ratio */
         //ASSERT(ytry >= *ylocmin);
         poly1->cursor = -1;
         while (structListReadNext(poly1, &p2, sizeof(FPOINT)) && p2.x < *xlocmin);
         // place the cursor in correspondence of the minimum
         READ_TRIPLET(poly1, p1, p2, p3, sizeof(FPOINT)); // updates p1, p2, p3
         // place the new point in the major interval:
         if (p2.x - p1.x > p3.x - p2.x)
            xtry = p2.x - GOLD*(p3.x - p2.x);
```

```c
            else
                xtry = p2.x + GOLD*(p2.x - p1.x);
            //ASSERT(inside(xtry, p1.x, p3.x));
            ytry = TO_BE_EVALUATED;
            newp = polyAddPoint(poly1, xtry, &ytry, 0, xtol);
            if (newp) {
                changes++;
                if (ytry < *ylocmin) {
                    *xlocmin = xtry, *ylocmin = ytry;
                    if (ytry < *ymin) {
                        *xmin = xtry, *ymin = ytry;
                        nFailed = 0;
                    }
                    else
                        nFailed++;
                    // N.B. the cursor of poly1 still points to the new point just added
                    READ_TRIPLET(poly1, p1, p2, p3, sizeof(FPOINT));
                }
                else { // place the cursor on current minimum, which is still p2!
                    nFailed++;
                    poly1->cursor = -1;
                    while (structListReadNext(poly1, &p2, sizeof(FPOINT)) && p2.x < *xlocmin);
                    READ_TRIPLET(poly1, p1, p2, p3, sizeof(FPOINT));
                }
                p2.marked = 1; structListWrite(poly1, &p2, sizeof(FPOINT));
                continue;
            }
        }
        else { /* !newp, point was not new,
                  if we already have 4 points we can try a cubic interpolation */
            if (poly1->len > 3) {
                double x02, x24;
                /* p1, p2, p3 are still valid, take as IV point p0 or p4,
                   the closest to p2, after p1 and p3 */
                poly1->cursor = -1;
                while (structListReadNext(poly1, &p2, sizeof(FPOINT)) && p2.x < *xlocmin);
                // cursor is placed on p2, still the current minimum
                x02 = structListRead(poly1, &p0, sizeof(FPOINT), -2) ? p2.x - p0.x : HUGE;
                // p0 is the point preceding p1 in the list
                x24 = structListRead(poly1, &p4, sizeof(FPOINT), 2) ? p4.x - p2.x : HUGE;
                // p4 is the point following p3 in the list
                if (x02 < x24) {
                    xtry = cubicMin(&p0, &p1, &p2, &p3, &ycast);
                    //ASSERT(inside(xtry, p0.x, p3.x));
                }
                else {
                    xtry = cubicMin(&p1, &p2, &p3, &p4, &ycast);
                    //ASSERT(inside(xtry, p1.x, p4.x));
                }
                if (xtry >= HUGE)
                    goto gold2; // cubic interpolation does not provide a minimum
            }
            else { // subdivide interval in golden ratio as a last resort
gold2:
                poly1->cursor = -1;
                while (structListReadNext(poly1, &p2, sizeof(FPOINT)) && p2.x < *xlocmin);
                // cursor is placed on p2, still the current minimum
                READ_TRIPLET(poly1, p1, p2, p3, sizeof(FPOINT));
                // place the new point in the major interval:
                if (p2.x - p1.x > p3.x - p2.x)
                    xtry = p2.x - GOLD*(p3.x - p2.x);
                else
                    xtry = p2.x + GOLD*(p2.x - p1.x);
                //ASSERT(inside(xtry, p1.x, p3.x));
            }
            ytry = TO_BE_EVALUATED;
            newp = polyAddPoint(poly1, xtry, &ytry, 0, xtol);
            if (newp) {
                changes++;
                if (ytry < *ylocmin) {
                    *xlocmin = xtry, *ylocmin = ytry;
                    if (ytry < *ymin) {
                        *xmin = xtry, *ymin = ytry;
                        nFailed = 0;
```

```c
            }
            else
                nFailed++;
            // N.B. the cursor of poly1 still points to the new point just added
            READ_TRIPLET(poly1, p1, p2, p3, sizeof(FPOINT));
        }
        else { // place the cursor on current minimum, which is still p2!
            nFailed++;
            poly1->cursor = -1;
            while (structListReadNext(poly1, &p2, sizeof(FPOINT)) && p2.x < *xlocmin);
            READ_TRIPLET(poly1, p1, p2, p3, sizeof(FPOINT));
        }
        p2.marked = 1; structListWrite(poly1, &p2, sizeof(FPOINT));
    }
  } while (newp /* loop as long as new points are added to the polygonal,
                   otherwise we have a static situation, i.e. we achieved convergence */
           && nFailed < nMaxFailed // otherwise we stop!
  );

  /* if the minimum is not improved (i.e. it's still the input guess)
     we must ensure that it's marked anyway (not to re-examine it) */
  poly1->cursor = -1;
  while (structListReadNext(poly1, &p2, sizeof(FPOINT)) && p2.x < *xlocmin);
  structListRead(poly1, &p2, sizeof(FPOINT), 0);
  p2.marked = 1; structListWrite(poly1, &p2, sizeof(FPOINT));

  return changes;
}

//---------------------------------------------------------------------
// Adds a point to the current polygonal approximation of the function:

int
polyAddPoint(STRUCTLIST *poly, double x, double *y, int m, double xtol)
/* If the point (x, y) is distinct from already present points
   the routine adds the point to the polygonal and returns 1,
   otherwise it does nothing and returns 0 */
{
  FPOINT p;
  // find out were to place the point in the right order:
  poly->cursor = -1;
  while (structListReadNext(poly, &p, sizeof(FPOINT))) {
    if (reqr(x, p.x, xtol)) {
      return 0; // a point with abscissa ~x is already present
      }
    if (p.x > x)
      break; // cursor is positioned on the point which must be next the new point
  }
  p.x = x;
  p.y = (*y == TO_BE_EVALUATED ? *y = (*f1)(x) // the function is evaluated HERE!
                               : *y);
  p.marked = m;
  structListInsertElement(poly, &p, sizeof(FPOINT));
  /* insert the new point before the element pointed by the cursor,
     this is left unchanged, so that it points to the freshly inserted element */
  return 1;
}

//---------------------------------------------------------------------
// Adjust a point in proximity of the domain boundary:

double
adjustToBoundary(double xa, double xb, double xc, double minRatio, int *clip)
/* Adjust xc with respect to the domain [xa, xc] where the function is defined:
   if xc is outside the domain it is returned the closest limit value and clip=1 is set,
   if xc is inside the domain, but close to a limit,
         this limit value is returned, but clip=0,
   otherwise the same value of xc is returned and clip=0 */
{
  if (xb > xa) {
    if (xc >= xsup) {
      *clip = 1;
      return xsup; // xc = xsup;
```

```c
      }
      else if ((xsup - xc)/(xc - xb) < minRatio) {
        *clip = 0;
        return xsup; // xc = xsup;
      }
      *clip = 0;
      return xc;
    }
    else {
      if (xc <= xinf) {
        *clip = 1;
        return xinf; // xc = xinf;
      }
      else if ((xinf - xc)/(xc - xb) < minRatio) {
        *clip = 0;
        return xinf; // xc = xinf;
      }
      *clip = 0;
      return xc;
    }
}

//--------------------------------------------------------------------------
// Computes the minimum of a parabola through three points:

double
quadricMin(FPOINT *p1, FPOINT *p2, FPOINT *p3, double *P2ymin)
/* Computes the point of minimum of a parabola passing through points p1, p2, p3,
   the abscissa is returned, and P2ymin is set to the ordinate.
   The Lagrange interpolation formula is used. */
{
#define xa (p1->x)
#define fa (p1->y)
#define xb (p2->x)
#define fb (p2->y)
#define xc (p3->x)
#define fc (p3->y)

#define A(fa,xa,xb,xc) (fa/(xa-xb)/(xa-xc))
#define B(fb,xa,xb,xc) (fb/(xb-xa)/(xb-xc))
#define C(fc,xa,xb,xc) (fc/(xc-xa)/(xc-xb))

#define P2(x,xa,fa,xb,fb,xc,fc)   A(fa,xa,xb,xc)*(x-xb)*(x-xc) \
                                + B(fb,xa,xb,xc)*(x-xa)*(x-xc) \
                                + C(fc,xa,xb,xc)*(x-xa)*(x-xb)
    double P2xmin;

    if (A(fa,xa,xb,xc) + B(fb,xa,xb,xc) + C(fc,xa,xb,xc) > 0.) // condition for a minimum
        P2xmin = (A(fa,xa,xb,xc)*(xb+xc)+B(fb,xa,xb,xc)*(xa+xc)+C(fc,xa,xb,xc)*(xa+xb))/2/(A(fa,xa,xb,xc)
+B(fb,xa,xb,xc)+C(fc,xa,xb,xc));
    else // we have a maximum instead
        P2xmin = HUGE; // a conventional very large and unlikely value

    *P2ymin = P2(P2xmin,xa,fa,xb,fb,xc,fc);

    return P2xmin;
}

//--------------------------------------------------------------------------
// Computes the minimum of a cubic through four points:

double
cubicMin(FPOINT *p1, FPOINT *p2, FPOINT *p3, FPOINT *p4, double *P3ymin)
/* Computes the point of minimum of a cubic passing through points p1, p2, p3, p4,
   the abscissa is returned, and P3ymin is set to the ordinate.
   The Lagrange interpolation formula is used. */
{
#define xa (p1->x)
#define fa (p1->y)
#define xb (p2->x)
#define fb (p2->y)
#define xc (p3->x)
#define fc (p3->y)
#define xd (p4->x)
```

```c
#define fd (p4->y)

#define AA(fa,xa,xb,xc,xd)  (fa/(xa-xb)/(xa-xc)/(xa-xd))
#define BB(fb,xa,xb,xc,xd)  (fb/(xb-xa)/(xb-xc)/(xb-xd))
#define CC(fc,xa,xb,xc,xd)  (fc/(xc-xa)/(xc-xb)/(xc-xd))
#define DD(fd,xa,xb,xc,xd)  (fd/(xd-xa)/(xd-xb)/(xd-xc))
#define A_(xa,fa,xb,fb,xc,fc,xd,fd)   (3*(AA(fa,xa,xb,xc,xd)+BB(fb,xa,xb,xc,xd)+CC(fc,xa,xb,xc,xd) \
+DD(fd,xa,xb,xc,xd)))
#define B_(xa,fa,xb,fb,xc,fc,xd,fd)   (AA(fa,xa,xb,xc,xd)*(xb+xc+xd)+BB(fb,xa,xb,xc,xd)*(xa+xc+xd) \
+CC(fc,xa,xb,xc,xd)*(xa+xb+xd)+DD(fd,xa,xb,xc,xd)*(xa+xb+xc))
#define C_(xa,fa,xb,fb,xc,fc,xd,fd)   (AA(fa,xa,xb,xc,xd)*(xb*xc+xb*xd+xc*xd) \
+BB(fb,xa,xb,xc,xd)*(xa*xc+xa*xd+xc*xd)+CC(fc,xa,xb,xc,xd)*(xa*xb+xa*xd+xb*xd) \
+DD(fd,xa,xb,xc,xd)*(xa*xb+xa*xc+xb*xc))
#define QQ(xa,fa,xb,fb,xc,fc,xd,fd) (SQR_(aaa, B_(xa,fa,xb,fb,xc,fc,xd,fd))- \
A_(xa,fa,xb,fb,xc,fc,xd,fd)*C_(xa,fa,xb,fb,xc,fc,xd,fd))
#define u31(xa,fa,xb,fb,xc,fc,xd,fd)   ((B_(xa,fa,xb,fb,xc,fc,xd,fd) \
+sqrt(QQ(xa,fa,xb,fb,xc,fc,xd,fd)))/A_(xa,fa,xb,fb,xc,fc,xd,fd))
#define u32(xa,fa,xb,fb,xc,fc,xd,fd)   ((B_(xa,fa,xb,fb,xc,fc,xd,fd)- \
sqrt(QQ(xa,fa,xb,fb,xc,fc,xd,fd)))/A_(xa,fa,xb,fb,xc,fc,xd,fd))
#define P3ii(x,fa,fb,fc,fd,xa,xb,xc,xd)   (A_(xa,fa,xb,fb,xc,fc,xd,fd)*x-B_(xa,fa,xb,fb,xc,fc,xd,fd))

#define P3(x,xa,fa,xb,fb,xc,fc,xd,fd) AA(fa,xa,xb,xc,xd)*(x-xb)*(x-xc)*(x-xd) \
                                    + BB(fb,xa,xb,xc,xd)*(x-xa)*(x-xc)*(x-xd) \
                                    + CC(fc,xa,xb,xc,xd)*(x-xa)*(x-xb)*(x-xd) \
                                    + DD(fd,xa,xb,xc,xd)*(x-xa)*(x-xb)*(x-xc)

    double P3xmin;

    if (QQ(xa,fa,xb,fb,xc,fc,xd,fd) > 0.)
        P3xmin = P3ii(u31(xa,fa,xb,fb,xc,fc,xd,fd),fa,fb,fc,fd,xa,xb,xc,xd) > 0. ?
                    u31(xa,fa,xb,fb,xc,fc,xd,fd) : u32(xa,fa,xb,fb,xc,fc,xd,fd); // conditions for a minimum
     else // we have a maximum instead
        P3xmin = HUGE; // a conventional very large and unlikely value

    *P3ymin = P3(P3xmin,xa,fa,xb,fb,xc,fc,xd,fd);

    return P3xmin;
}

//---------------------------------------------------------------------------
// Checks if two floating point numbers are almost equal in relative terms:

int
reqr(double a, double b, double eps)
{
    return fabs(a - b) < eps*(1 + fabs(a));
}

//---------------------------------------------------------------------------
// Evaluates the function at x and updates the number of function evaluations:

static double
f1(double x)
// Evaluates the function at x and increments the number of function evaluations
{
    double fx;   //ASSERT(contained(x, xinf, xsup));
    ++nff; fx = (*ff)(x);
    return fx;
}

//---------------------------------------------------------------------------
// test if x is strictly contained in [a, b]
double inside(double x, double a, double b)
{
   return x > a && x < b;
}

//---------------------------------------------------------------------------
// no operation:
void
noop()
{}
```

# 6    Conclusions

We proposed a new algorithm for the optimization of one-dimensional, arbitrary functions. The theory of operation, the results obtained in some practical cases and the C code of the algorithm have been illustrated and discussed in detail. We think that the algorithm improves over the classical techniques described in literature in many aspects, it also takes into consideration a wide range of possible complications that may arise, thus it should be suitable for general use in almost any practical case.

(*) Glauco Masotti was born in Ravenna, Italy, in 1955. He graduated summa cum laude Doctor in Electronic Engineering at the University of Bologna, Italy, in 1980.
Since then he has been involved in the design and development of various software systems, but he also got patents for some electronic devices.
He specialized in CAD/CAM systems, and in research in the field. Most of his work, being proprietary, has not been published. In early years, in Bologna, he cooperated at the development of the GBG drafting system for CAD.LAB (former name of Think3).
For COPIMAC he led the development of Aliseo, one of the first highly interactive, parametric and feature-based, solid modelers.
In 1989 he was visiting scholar at the University of Southern California (USC), in Los Angeles, where he worked at the design of an object-oriented geometric toolkit. He also proposed an extension to the C++ language, known as EC++ (Extended C++).
Back in Italy he undertook research at the Dept. of Electronics and Computer Science of the University of Bologna on symbolic methods and numerical problems in geometric algorithms. He also taught a course on Industrial Automation, with emphasis on geometric modeling and CAD/CAM integration.
In 1992 he joined Ecocad (later acquired by Think3), in Pesaro, Italy, where he was in charge of the design and development of various parts of Eureka, an advanced surface and solid modeler. One of his major contributions was the development of a module for assembly modeling and kinematic simulation, supporting the positioning of parts via mating constraints. In 1997 he was co-founder of IeS, Ingegneria e Software Srl, in Bologna, where he worked till 2001, developing software for automatic conversion of drawings in solid models, automatic production of exploded views, interpolation of measured 3D points, plane development of surfaces with minimal distortion, and finite element analysis. He was also involved in structural analysis and CFD problems.
In latest years, while trading for a living, he undertook autonomous researches in various fields, mainly on mathematical optimization and digital filters.
He is an enthusiast of sailing, which he practiced at a competitive level until 2002, racing on catamarans.